\documentclass[11pt]{amsart}
\usepackage{amsthm,amsmath,epsf,amssymb,enumerate,array}
\usepackage[all]{xy} 

\newtheorem{thm}{Theorem}[section]
\newtheorem{lemma}[thm]{Lemma}

\newtheorem{claim}[thm]{Claim}

\newtheorem{prop}[thm]{Proposition}

\newtheorem{remark}[thm]{Remark}
\theoremstyle{plain}

\DeclareFontFamily{U}{rsf}{}
\DeclareFontShape{U}{rsf}{m}{n}{<5><6>rsfs5<7><8><9> rsfs7 <10-> rsfs10}{}
\DeclareMathAlphabet\Scr{U}{rsf}{m}{n}
\DeclareMathAlphabet\mathbi{U}{cmr}{bx}{it}

\def\CY{Calabi-Yau}
\def\roof{\mbox{\tiny \mbox{$\!\vee$}}}
\def\comp{\mbox{\scriptsize \mbox{$\,\circ\,$}}}
\def\c#1{\mathcal{#1}}

\def\P{{\mathbb P}}

\def\Z{{\mathbb Z}}

\def\D{\mathrm{D}}

\def\Ll{\mathbf{L}}

\def\iso{\cong}
\def\niso{\ncong}
\def\bar{\overline}

\def\too{\longrightarrow}

\def\id{\operatorname{id}}
\def\Hom{\operatorname{Hom}}

\def\RHom{\R\!\operatorname{Hom}}

\def\Aut{\operatorname{Aut}}

\def\im{\operatorname{im}}

\def\dim{\operatorname{dim}}

\def\Ltensor{\mathbin{\overset{\mbox{\tiny\mbox{$\Ll$}}}\otimes}}

\def\ms#1{\mathsf{#1}}

\def\ses#1#2#3{\xymatrix@1{0 \ar[r] & #1 \ar[r] & #2 \ar[r] & #3 \ar[r] & 0}}

\DeclareMathOperator{\spec}{Spec} 

\newcommand{\mono}{\hookrightarrow} 
\newcommand{\mor}[1]{\xrightarrow{#1}} 
\newcommand{\isomor}{\mor{\sim}} 
\newcommand{\card}[1]{\lvert #1\rvert} 
\newcommand{\abs}[1]{\lvert#1\rvert} 
\newcommand{\dual}{^{\roof}} 
\newcommand{\mrs}{^{\sharp}} 
\newcommand{\sww}{\w}
\newcommand{\rest}[1]{\lvert_{#1}} 
\newcommand{\pr}{\pi} 
\newcommand{\diag}{\delta} 

\newcommand{\K}{\Bbbk} 
\newcommand{\cat}[1]{{\mathfrak{#1}}} 
\newcommand{\Coh}{\cat{Coh}} 
\newcommand{\s}[1]{\mathcal{#1}} 
\newcommand{\p}{{\rm P}} 
\newcommand{\w}{{\rm w}} 
\newcommand{\lotimes}{\Ltensor} 
\newcommand{\lboxtimes}{\mathbin{\overset{\mbox{\tiny\mbox{$\Ll$}}}\boxtimes}}
\newcommand{\lboxtimesa}{\mathbin{
\overset{\mbox{\tiny\mbox{$\Ll$}}}\boxtimes\,}}
\newcommand{\so}{\s{O}} 
\newcommand{\sd}{\Omega} 
\newcommand{\diff}[1]{d_{#1}} 
\newcommand{\sko}{\s K} 
\newcommand{\ind}[1]{_{#1}} 
\newcommand{\bo}[1]{\s{M}\ind{#1}} 
\newcommand{\abo}[1]{\tilde{\s{M}}\ind{#1}} 
\newcommand{\cone}[1]{\mathsf{C}\! \left( #1 \right)} 
\newcommand{\coneI}[1]{\mathsf{C} \left( #1 \right)} 
\newcommand{\mc}[1]{\mathsf{MC}\! \left( #1 \right)} 
\newcommand{\mcI}[1]{\mathsf{MC} \left( #1 \right)} 
\newcommand{\FM}[2][]{\Phi^{#1}_{#2}} 
\newcommand{\FMcomp}{\star\,} 
\newcommand{\dr}{\s{R}} 
\newcommand{\pdr}[1]{\dr\ind{#1}} 
\newcommand{\bdr}{\s{B}} 
\newcommand{\bpdr}[1]{\bdr\ind{#1}} 
\newcommand{\mpdr}[1]{\alpha\ind{#1}} 
\newcommand{\ndr}{\nu} 
\newcommand{\drim}[1]{\s{C}\ind{#1}} 
\newcommand{\mdrim}[1]{\ndr\ind{#1}} 
\newcommand{\G}{\s{G}} 
\newcommand{\Gi}[1]{\G\ind{#1}} 
\newcommand{\cprest}[1]{\widehat{#1}} 
\newcommand{\drX}{\cprest{\dr}} 
\newcommand{\pdrX}[1]{\drX\ind{#1}} 
\newcommand{\mpdrX}[1]{\cprest{\alpha}\ind{#1}} 
\newcommand{\Gim}[1]{\s{F}\ind{#1}} 
\newcommand{\GimX}[1]{\cprest{\s{F}}\ind{#1}} 
\newcommand{\mGi}[1]{\eta\ind{#1}} 
\newcommand{\amGi}[1]{\tilde{\eta}\ind{#1}} 
\newcommand{\boX}[1]{\cprest{\s{M}}\ind{#1}} 
\newcommand{\aboX}[1]{\cprest{\tilde{\s{M}}}\ind{#1}} 
\newcommand{\mbo}[1]{\epsilon\ind{#1}} 
\newcommand{\mboX}[1]{\cprest{\epsilon}\ind{#1}} 
\newcommand{\Times}{\!\times\!}
\newcommand{\Timesa}{\!\times}
\newcommand{\R}{{\mathbf R}} 
\newcommand{\KK}{\mathrm{K}} 
\newcommand{\CC}{\mathrm{C}} 
\newcommand{\lradj}[1]{\beta_{#1}} 
\newcommand{\rladj}[1]{\gamma_{#1}} 

\begin{document}

\title[Derived category autoequivalences and Beilinson]{
Derived autoequivalences and a \\ weighted Beilinson resolution}

\author[A. Canonaco]{Alberto Canonaco}
\address{Dipartimento di Matematica ``F. Casorati'', Universit{\`a}
di Pavia, Via Ferrata 1, 27100 Pavia, Italy}
\email{alberto.canonaco@unipv.it}

\author[R.L. Karp]{Robert L. Karp}
\address{NHETC, Rutgers University, 126 Frelinghuysen Rd, Piscataway,
NJ 08854-8019 USA}
\email{karp@rci.rutgers.edu}



\begin{abstract}
Given a smooth stacky Calabi-Yau hypersurface $X$ in a weighted
projective space, we consider the functor $\ms{G}$ which is the
composition of the following two autoequivalences of $\D^b(X)$: the
first one is induced by the spherical object $\so_X$, 
while the second one is
tensoring with $\so_X(1)$. The main result of the paper is that the
composition of $\ms{G}$ with itself $\sww$ times, where $\sww$ is the
sum of the weights of the weighted projective space, is isomorphic to
the autoequivalence ``shift by 2''. The proof also involves the
construction of a Beilinson type resolution of the diagonal for
weighted projective spaces, viewed as smooth stacks.
\end{abstract}

\maketitle
\vskip 1cm

\section{Introduction}    \label{s:intro}

Recent years have seen an increased activity in the study of algebraic
varieties via their derived categories of coherent sheaves. Although this
algebraic approach is indirect compared to the geometric investigations
involving divisors, curves, or branched covers, just to name a few, it is
nevertheless quite promising, as in some cases it allows for a deeper
understanding. This is the case of varieties with  interesting groups of
derived autoequivalences, and in particular those with Kodaira dimension 0,
where the autoequivalences are symmetries of the variety not visible in the
geometric presentation.

Despite significant progress, our understanding of the derived categories of
coherent sheaves and their autoequivalences is limited (for a recent review we
refer to \cite{Bridgeland:Review}). In the present paper we hope to further
this  understanding by proving certain identities involving Fourier-Mukai
functors on quasi-smooth \CY\ varieties, viewed as smooth Deligne-Mumford
stacks.
The origin of these identities is closely tied with mirror symmetry.
We first
present our main result, then describe the setting in which it arises.

Let $\P^n({\bf w})$ be an $n$-dimensional weighted projective space
over a field $\K$, regarded
as a smooth proper Deligne-Mumford stack, with weight vector ${\bf
w}=(\w_0,\w_1,\ldots,\w_n)$, and let $\sww = \sum_{i=0}^n \w_i$ denote the sum
of all the weights. We have several equivalent ways to think about $\P^n({\bf
w})$: graded scheme \cite{Alberto}, toric stack \cite{BorisovDM}, or quotient
stack \cite{Auroux}.

Let $X$ be an anti-canonical hypersurface in $\P^n({\bf w})$.
By the stacky version of Bertini's theorem
the generic member of the linear system $\abs{-K_{\P^n({\bf w})}}=
\abs{\so_{\P^n({\bf w})}(\sww)}$ is a  proper smooth Deligne-Mumford stack. Let
$\D^b(X)$ denote the bounded derived category of coherent sheaves on the stack
$X$. We have two functors naturally associated to this data:
\begin{equation}\label{k2}
   \ms{L}\colon \D^b(X)\to \D^b(X)\qquad \ms{L}(\c F) = \c F\otimes \so_X(1),
\end{equation}
and
\begin{equation}\label{k3}
   \ms{K}\colon \D^b(X)\to \D^b(X)\qquad \ms{K}(\c F) =
     \cone{\RHom_{X}(\so_X,\c F)\otimes_{\K}\so_X\too\c F},
\end{equation}
for all $\c F \in \D^b(X)$.
The morphism in $\ms{K}$ is the evaluation map, and $\mathsf{C}$ is its
cone. For an integer $m$ we also have the autoequivalence of $\D^b(X)$
given by the translation functor $(-)[m]$; its action is ``shift by
$m$'', i.e., $\c F\mapsto \c F[m]$. In fact also
$\ms{L}$ is clearly an equivalence, and the same is true for
$\ms{K}$ thanks to \cite{ST:braid}, where it is proved more generally
that the functor defined by $\c F\mapsto
\cone{\RHom_{X}(\c E,\c F)\otimes_{\K}\c E\too\c F}$ is an
autoequivalence of $\D^b(X)$ whenever $\c E$ is a {\em spherical
object}, i.e., an object of $\D^b(X)$ such that
$\c E\otimes\omega_X\iso\c E$ and
\[
\Hom_{\D^b(X)}(\c E,\c E[i])\iso\begin{cases}
\K & \text{if $i=0,\dim(X)$}\\
0 & \text{otherwise.}
\end{cases}
\]
Note that $X$ being \CY\ immediately implies that $\so_X$ is a
spherical object. Our main result is then the following non-trivial
relation in the group $\Aut(\D^b(X))$ of (isomorphism classes of)
autoequivalences of $\D^b(X)$:
\begin{thm}\label{t1}
Let $X$ be a smooth anti-canonical stacky hypersurface in the Deligne-Mumford
stack $\P^n({\bf w})$, and let $\ms{G}=\ms{L}\comp\ms{K}$, where
$\ms{L}$ and $\ms{K}$ are the autoequivalences of $\D^b(X)$ defined in
\eqref{k2} and \eqref{k3}. Then there is an isomorphism of functors:
\begin{equation}\label{eG}
\underbrace{\ms{G} \comp \cdots \comp \ms{G}}_{\txt{$\sww$-times}} \iso (-)[2].
\end{equation}
\end{thm}

Let us spend some time trying to understand the origin of this
statement. 
At first sight it might seem surprising that physics has anything to do with
such an abstract branch of pure mathematics. But one should remember that
historically some of the most interesting mathematical problems came from the
real world, and primarily from physics. Recently, with the advent of string
theory, the bridge of interactions between abstract mathematics and theoretical
physics has entered an era of renaissance, with mirror symmetry the most
prominent example of the interaction.

{From} the point of view of {\em strings} in string theory the
appearance of the derived category is quite intriguing, but recent
developments showed that D-branes mandate a categorical approach. In
particular, Douglas argued that B-type topological D-branes are
objects in the bounded derived category of coherent sheaves
\cite{Douglas:2000gi}. His work was subsequently axiomatized by
Bridgeland \cite{Bridgeland:2002}, and has since been subject to
active investigations.  The A-type D-branes have a very different
description, involving the derived Fukaya category. Mirror symmetry
exchanges the A and B branes, and naturally leads to Kontsevich's {\em
homological mirror symmetry} (HMS) conjecture.  For a detailed
exposition of these ideas we refer the reader to recent book
\cite{DBook}.

To motivate our result we need to start with mirror symmetry in its pre-HMS
phase. In this form mirror symmetry is an isomorphism between the
(complexified) Kahler moduli space $\c M_K(X)$ of a \CY\ variety $X$ and the
moduli space of complex deformations $\c M_c(\widetilde{X})$ of its mirror
$\widetilde{X}$.
For the precise definitions we refer to the book by Cox and
Katz \cite{Cox:Katz}. We will follow their terminology in this introduction. 

We note at this point that the moduli spaces in question are only coarse
moduli spaces; the fine moduli spaces are necessarily stacks.
This fact complicates any existing intuitive mathematical understanding, 
but the conformal field theory (CFT) techniques that underlie the Cox and
Katz exposition give us an alternative view, which we now elaborate on.

Mirror symmetry also suggested a natural way to complexify the Kahler
moduli space and how to compactify it. The complexified Kahler moduli
space $\c M_K(X)$ in general is an intricate object, but for $X$ a
hypersurface in a toric variety it has a rich combinatorial structure
and is relatively well-understood.  In particular, the fundamental
group of $\c M_K(X)$ in general is non-trivial, and one can talk about
various monodromy representations. More concretely, there are two
types of boundary divisors in $\c M_K(X)$: ``large radius divisor''
and the ``discriminant'' (some authors refer to both as discriminant,
but for us the distinction is important). Both of these are reducible
in general.  At the large radius divisor certain cycles of $X$, viewed
as a Kahler manifold, acquire infinite volume. The discriminant is
somewhat harder to describe.  The original definition is that the CFT
associated to a string probing $X$ becomes singular at such a point in
moduli space.  Generically this happens because some D-brane (or
several of them, even infinitely many) becomes massless, and therefore
the effective CFT description provided by the string fails. A
consequence of this fact is that, by using the mirror map isomorphism
of the moduli spaces, as one approaches the discriminant in $\c
M_K(X)$ one is moving in $\c M_c(\widetilde{X})$ to a point where the
mirror $\widetilde{X}$ is developing a singularity.

Armed with this picture of $\c M_K(X)$, we can fix a basepoint $O$,
and look at loops in $\c M_K(X)$ based at $O$.  The CFT description of
string theory shows that traversing such a loop gives in general a
non-trivial functor $\D^b(X)\to \D^b(X)$, which moreover has to be an
equivalence (string theory does not seem to able to distinguish
between isomorphism and equivalence).  Therefore we arrive at a group
homomorphism, first suggested by Kontsevich \cite{Kont:mon}:
\begin{equation*}
 \mu \colon \pi_1(\c M_K(X)) \too \Aut(\D^b(X)).
\end{equation*}
At present writing very little is known about $\mu$. Kontsevich's
ideas were generalized by Horja \cite{Horj:DX} and Morrison
\cite{Mor:geom2}.
The question at hand is: given a pointed loop
in $\c M_K(X)$, what is the associated autoequivalence in $ \D^b(X)$? 
Progress in this direction was made in \cite{en:Horja}, 
where this question is answered for the EZ-degenerations 
introduced in \cite{Horj:EZ}. 

It is clear now that 
given a presentation of $\pi_1(\c M_K(X))$ 
where we know the images under $\mu$ of the generators, 
the relations in the presentation will determine interesting identities 
in $\Aut(\D^b(X))$. We now turn to an example of this sort. 
Let $X$ be a smooth
degree $\sww =\sum_{i=0}^n\w_i$ variety in $\P^n({\bf w})$, in other words
let $X$ be such that it
does not meet the singularities of $\P^n({\bf w})$. In this case the
compactification of $\c M_K(X)$ is isomorphic to $\P^1$, and we have
three distinguished
points $P_{LC}$, $P_0$ and $P_{F}$. It is easier to describe them in terms of
$\c M_c(\widetilde{X})$: $P_{LC}$ is a large complex structure limit point
(with maximally unipotent monodromy), at $P_0$ the family
$\widetilde{X}$ has rational double points,  while at $P_{F}$ it has
additional automorphisms. If $\widetilde{X}$ has a Fermat form, i.e.,
$\w_i$ divides $\sww$, then we are talking about
an additional cyclic symmetry $\Z_\sww$.

Let $M_P$ denote the monodromy associated to a loop around the point
$P$. Since $P_{LC}$ and $P_0$ are the only limit points of $\c
M_K(X)$, and the compactification of this is isomorphic to $\P^1$ (see
\cite{Cox:Katz}), with $\pi_1(\P^1- \{\mbox{2 points} \})=\Z$, one
would want to conclude, incorrectly, that $M_{P_{LC}}$ and $M_{P_0}$
are related.  On other hand, the extra automorphisms indicates that
$P_{F}$ is a stacky point in the moduli space, with finite stabilizer,
and so, at best, the $\sww$-th power of $M_{P_{F}}\iso M_{P_{LC}}\comp
M_{P_0}$ is the identity.  In the case of $\P^n({\bf w})=\P^4$
Kontsevich proposed that $M_{P_0}=\ms{K}$ from \eqref{k3},
$M_{P_{LC}}=\ms{L}$, and checked that indeed $M_{P_{F}}^5=\id$ in
K-theory.  Later on Aspinwall realized that in fact $M_{P_{F}}^5\iso
(-)[2]$. Based on physical considerations it was clear to us that the
Kontsevich-Aspinwall result should hold in the weighted case as well,
which eventually led us to Theorem~\ref{t1}. Cases where $\c M_K(X)$
is higher dimensional were investigated in \cite{K1,K2}.


Our proof is inspired by \cite[Sec.~7.1.4]{Paul:TASI2003}, where the
case $\P^n({\bf w})=\P^4$ is outlined. Actually with the same
technique a more general result was independently obtained by
Kuznetsov in \cite[\S 4]{V14}, where smooth Fano (but the argument
applies to the Calabi-Yau case as well) hypersurfaces in $\P^n$ were
considered (in Remark~\ref{Fano} we explain how Kuznetsov's result
extends to the weighted case, as suggested to us by the author after
the first version of this paper was made public).  However, the proof
of Theorem~\ref{t1} is much harder, and a different approach is needed
overall.  Still, the idea to trade the composition of the functors for
the composition of their kernels, and then use a resolution of the
diagonal of the (weighted in our case) projective space proved very
useful to us. In fact a good part of the paper is devoted to the
construction of a resolution of the diagonal for $\P^n({\bf w})$,
which is similar but still very different from the well known
Beilinson resolution of $\P^n$.

\subsubsection*{\bf Acknowledgments}
It is a pleasure to thank Tom Bridgeland, Mike Douglas, Sheldon Katz,
Alastair King, Alexander Kuznetsov, David Morrison, Tony
Pantev,  Ronen Plesser and especially Paul Aspinwall for useful conversations.
We would also like to thank the 2005 Summer Institute in Algebraic Geometry at
the University of Washington, for providing a stimulating environment where
this work was initiated.

\subsubsection*{\bf Notation}
A complex $A$ of some abelian category $\cat{A}$ is given by a collection of
objects $A^i$ of $\cat{A}$, together with morphisms $\diff{A}^i\colon A^i\to
A^{i+1}$ such that $\diff{A}^{i+1}\comp\diff{A}^i=0$ for every $i\in\Z$. When
$A$ is just an object of $\cat{A}$, it is viewed as a bounded complex with
$A^0=A$ and $A^i=0$ for $i\ne0$. A morphism of complexes $f\colon A\to B$ is
given by a collection of maps $f^i\colon A^i\to B^i$ such that $\diff{B}^i\comp
f^i=f^{i+1}\comp\diff{A}^i$. For $k\in\Z$ the shifted complex $A[k]$ is defined
by $A[k]^i:=A^{i+k}$ and $\diff{A[k]}^i:=(-1)^k\diff{A}^{i+k}$; similarly,
$f[k]$ is defined by $f[k]^i:=f^{i+k}$.

The dual of a locally free sheaf $\s L$ (or more generally of a complex of
locally free sheaves) will be denoted by $\s L\dual$; the same notation will be
used for the dual of a vector space.

\tableofcontents

\section{Preliminaries on Fourier-Mukai functors}    \label{s:1}

In order to fix notations and conventions that will be used throughout
the paper, we start recalling some definitions and basic facts about
triangulated categories, derived categories and derived functors. For
a thorough treatment of the subject we refer to \cite{Hart:dC} or
\cite{Huybrechts}.

If $f\colon A\to B$ is a morphism in a triangulated category, then a {\em cone}
of $f$ is an object $\cone{f}$ (defined up to isomorphism), which fits into a
distinguished triangle
\[\xymatrix@1{A\ar[r]^-{f} & B \ar[r] & \cone{f} \ar[r] & A[1]}.\]

For an abelian category $\cat{A}$, $\CC^b(\cat{A})$ denotes the abelian
category of bounded complexes of $\cat{A}$ (its objects are complexes $A$ such
that $A^i=0$ for $\abs{i}\gg0$). The {\em mapping cone} of a morphism $f\colon
A\to B$ of $\CC^b(\cat{A})$ is the complex $\mc{f}$ defined by
\begin{equation}\label{vv1}
\mc{f}^i:=A^{i+1}\oplus B^i\,,\qquad
\diff{\mcI{f}}^i:=\begin{pmatrix}
-\diff{A}^{i+1} & 0\\
f^{i+1} & \diff{B}^i
\end{pmatrix}.
\end{equation}
$\KK^b(\cat{A})$ (respectively $\D^b(\cat{A})$) will be the bounded homotopy
(respectively derived) category of $\cat{A}$. We recall that $\KK^b(\cat{A})$
and $\D^b(\cat{A})$ are triangulated categories and have the same objects as
$\CC^b(\cat{A})$. For a morphism $f$ in $\CC^b(\cat{A})$, $f$ will
also denote its image in $\KK^b(\cat{A})$, or in $\D^b(\cat{A})$. In both
categories $\cone{f}\iso\mc{f}$, but $\mc{f}$ will be used only when the
specific form of the resulting complex is needed.

If $\cat{B}$ is another abelian category, a left exact functor
$F\colon\cat{A}\to\cat{B}$ trivially extends to an exact functor again denoted
by $F\colon\KK^b(\cat{A})\to\KK^b(\cat{B})$. When it exists, its right derived
functor will be denoted by $\R F\colon\D^b(\cat{A})\to\D^b(\cat{B})$; we set
$R^iF:=H^i\comp\R F$ for $i\in\Z$. If $A$ is an object of $\D^b(\cat{A})$ such
that each $A^i$ is $F$-acyclic (i.e., $R^jF(A^i)=0$ for $j>0$), then $\R
F(A)\iso F(A)$ in $\D^b(\cat{B})$. Similar considerations hold if $F$ is right
exact, in which case its left derived functor will be denoted by $\Ll F$.

For simplicity in the following we will call {\em stack} 
a Deligne-Mumford stack which is proper and smooth over the base field
$\K$, and such that every coherent sheaf is a quotient of a locally
free sheaf of finite rank. In fact all stacks we will consider in the
rest of the paper will be stacks associated to normal projective
varieties with only quotient singularities (namely, weighted
projective spaces, quasi-smooth hypersurfaces in them and products of
such varieties), and those satisfy our condition thanks to
\cite[Theorem 4.2]{Kawamata:DC}. When the proofs remain essentially
the same, we will use results stated in the literature only
for schemes for stacks as well, but most of the time we point this out.

If $Y$ is a stack, $\Coh(Y)$ will denote the abelian category of
coherent sheaves on $Y$, and we set for brevity
$\CC^b(Y):=\CC^b(\Coh(Y))$, $\KK^b(Y):=\KK^b(\Coh(Y))$ and
$\D^b(Y):=\D^b(\Coh(Y))$. If $f\colon Y\to Z$ is a morphism of stacks,
there are derived functors $\R f_*\colon\D^b(Y)\to\D^b(Z)$ and $\Ll
f^*\colon\D^b(Z)\to\D^b(Y)$. Notice that $\R f_*\iso f_*$ if $f$ is
finite and $\Ll f^*\iso f^*$ if $f$ is flat. When $Z$ is a point $f_*$
can be identified with $\Gamma(Y,-)$, and $R^i\Gamma(Y,-)$ will be
denoted by $H^i(Y,-)$. Our definition of stack also implies that there
is a left derived functor for the tensor product, denoted by
$-\lotimes-\colon\D^b(Y)\Times\D^b(Y)\to\D^b(Y)$. Clearly $\s
F\lotimes\s G\iso\s F\otimes\s G$ if, either, each $\s F^i$, or each
$\s G^i$, is locally free.

Given $\s{E}\in \D^b(Y\Times Z)$, and denoting by $p\colon Y\Times Z\to Y$ and
$q\colon Y\Times Z\to Z$ the projections, the exact functor $\FM{\s E}=\FM[Z\to
Y]{\s E}\colon\D^b(Z)\to\D^b(Y)$ defined by
\[\FM{\s E}(\s F):=\R p_*(\s E\lotimes q^*\s F)\] for $\s F\in \D^b(Z)$, is
called a {\em Fourier-Mukai functor} with kernel $\s{E}$.

\begin{lemma}\label{ker0}
If $\s E\in \D^b(Y\Times Z)$ is such that $\FM{\s E}\iso0$, then $\s E\iso0$.
\end{lemma}

\begin{proof}
Assume on the contrary that $\s E\niso0$, and let $m$ be the least integer such
that $H^m(\s E)\ne0$. Setting $\s F:=H^m(\s E)\in\Coh(Y\Times Z)$, we claim
that it is enough to prove that
\begin{equation}\label{eL}
p_*(\s F\otimes q^*\s L)\ne0\qquad\text{for some $\s L\in\Coh(Z)$ locally
free}.
\end{equation}
Indeed, assuming this, the definition of $\s F$ implies that there is
a distinguished triangle in $\D^b(Y\Times Z)$
\[\xymatrix@1{\s E'[-1] \ar[r] & \s F[-m] \ar[r] & \s E \ar[r] & \s E'}\]
with $H^i(\s E')=0$ for $i\le m$, from which it is easy to deduce that
$R^ip_*(\s E'\otimes q^*\s L)=0$ for $i\le m$. Then, applying the
exact functor $\R p_*(-\otimes q^*\s L)$ to the above triangle, and
taking the associated cohomology sequence, we obtain
\[
0\ne p_*(\s F\otimes q^*\s L)\iso R^mp_*(\s F[-m]\otimes q^*\s L)
\iso R^mp_*(\s E\otimes q^*\s L)\iso H^m(\FM{\s E}(\s L)),
\]
which contradicts the hypothesis $\FM{\s E}\iso0$.

In order to prove \eqref{eL} it is obviously enough to find a locally
free sheaf $\s L$ such that $f^*p_*(\s F\otimes q^*\s L)\ne0$, where
$f\colon V\to Y$ is an \'etale and surjective morphism and $V$ an
affine scheme. Applying the ``flat base change'' theorem (for stacks
this is \cite[Prop. 13.1.9]{Stacks}) to the Cartesian square
\[
\xymatrix{V\Times Z \ar[rr]^-{\bar f} \ar[d]_{\bar p} & & Y\Times Z \ar[d]^p \\
V \ar[rr]^-f & & Y}
\]
and setting $\bar{\s F}:=\bar f^*\s F\in\Coh(V\Times Z)$, we see that
$f^*p_*(\s F\otimes q^*\s L)\iso
\bar p_*(\bar{\s F}\otimes\bar q^*\s L)$, 
where $\bar q=q\comp\bar f\colon V\Times Z\to Z$ is the
projection. Hence it is enough to show that
\[
0\ne\Hom_V(\so_V,\bar p_*(\bar{\s F}\otimes\bar q^*\s L))\iso\Hom_{V\Timesa
Z}(\so_{V\Timesa Z},\bar{\s F}\otimes\bar q^*\s L)\iso\Hom_{V\Timesa Z}(\bar
q^*\s L\dual,\bar{\s F})\iso\Hom_Z(\s L\dual,\bar q_*\bar{\s F}).
\]
Taking into account that $\bar q_*\bar{\s F}=\bar q_*\bar f^*\s F\ne0$
(because $\bar q$ is affine and $\bar f$ is \'etale and surjective)
and the fact that the quasi-coherent sheaf $\bar q_*\bar{\s F}$ is the
inductive limit of its coherent subsheaves (by
\cite[Prop. 15.4]{Stacks}), the existence of $\s L\in\Coh(Z)$ locally
free such that $\Hom_Z(\s L\dual,\bar q_*\bar{\s F})\ne0$ follows from
the fact that every coherent sheaf on $Z$ is a quotient of a locally
free sheaf.
\end{proof}

Now we specialize to the case $Y=Z$, although much of what we are
going to say can be extended to the general case, with obvious
modifications.

Given $\s{F},\s{G}\in\KK^b(Y)$, their {\em exterior} tensor product is
defined as
\begin{equation*}
\s{F}\boxtimes\s{G}:=\pr_2^*{\s{F}}\otimes\pr_1^*{\s{G}}\in\KK^b(Y\Times Y),
\end{equation*}
where $\pr_i$ is the natural projection to the $i$th factor of $Y\Times Y$. The
symbol $\lboxtimes$ will be used for the derived functor of exterior tensor
product (again, $\s F\lboxtimes\s G\iso\s F\boxtimes\s G$ if either
each $\s F^i$ or each $\s G^i$ is locally free).

Given $\s{E}\in \D^b(Y\Times Y)$, we define two Fourier-Mukai functors
$\FM[i]{\s{E}}\colon\D^b(Y)\to\D^b(Y)$:
\begin{equation*}
\FM[i]{\s{E}}(\s{F}):=\R{\pr_{3-i}}_*(\s{E}\lotimes\pr_i^*\s{F})
,\qquad \text{for $i=1,2$}.
\end{equation*}

The composition of Fourier-Mukai functors is again a Fourier-Mukai
functor (see \cite[Prop. 5.10]{Huybrechts}). In
particular,
\begin{equation*}
\FM[1]{\s{F}\FMcomp\s{E}}\iso \FM[1]{\s{F}}\, \comp\, \FM[1]{\s{E}},
\end{equation*}
where $\FMcomp$ is the composition of kernels, and for $\s{E},\s{F}\in
\D^b(Y\Times Y)$ it is defined by
\begin{equation*}
\s{F}\FMcomp\s{E}=\R{\pr_{1,3}}_*(\pr_{2,3}^*\s{F}\lotimes\pr_{1,2}^*\s{E})
\in \D^b(Y\Times Y).
\end{equation*}
$\pr_{i,j}$ is projection to the $i$th times $j$th factor in the
product $Y\Times Y\Times Y$.

\begin{lemma}\label{kercomp}
Given $\s{F},\s{E}_i,\s{F}_i\in\D^b(Y)$, for $i=1,2$, and $\c W \in
\D^b(Y\Times Y)$, we have the following isomorphisms:
\begin{enumerate}
\item $\FM[i]{\s{E}_2\lboxtimesa\s{E}_1}(\s{F})\iso
\s{E}_{3-i}\otimes_{\K}\R\Gamma(Y,\s{E}_i\lotimes\s{F})$, for $i=1,2$;
\item $(\s{F}_2\lboxtimes\s{F}_1)\FMcomp\c W\iso
\s{F}_2\lboxtimes\FM[2]{\c W}(\s{F}_1)$;
\item $\c W\FMcomp (\s{F}_2\lboxtimes\s{F}_1)\iso
\FM[1]{\c W}(\s{F}_2)\lboxtimes\s{F}_1$;
\item $(\s{F}_2\lboxtimes\s{F}_1)\FMcomp(\s{E}_2\lboxtimes\s{E}_1)\iso
\s{F}_2\lboxtimes\s{E}_1\otimes_{\K}\R\Gamma(Y,\s{E}_2\lotimes\s{F}_1)$.
\end{enumerate}
\end{lemma}

\begin{proof}
(1) By definition we have
\begin{equation*}
\FM[i]{\s{E}_2\lboxtimesa\s{E}_1}(\s{F})\iso
\R{\pr_{3-i}}_*\left( \pr_{3-i}^*\c E_{3-i}\lotimes\pr_i^*(\c E_i
\lotimes{\s{F}})\right)
\end{equation*}
and using the projection formula
\begin{equation*}
\R{\pr_{3-i}}_*\left( \pr_{3-i}^*\c E_{3-i}\lotimes\pr_i^*(\c E_i
\lotimes{\s{F}})\right)\iso
\c E_{3-i}\lotimes\R{\pr_{3-i}}_*\pr_i^*(\c E_i \lotimes{\s{F}}).
\end{equation*}
Finally, $\R{\pr_{3-i}}_*\pr_i^*(\c E_i\lotimes\s{F}))\iso
\so_Y\otimes_{\K}\R\Gamma(Y,\s{E}_i\lotimes\s{F})$ by the ``flat base change''
theorem applied to the Cartesian square
\begin{equation}\label{e22}
\xymatrix{
Y\Times Y \ar[rr]^{\pr_1} \ar[d]_{\pr_2} & & Y \ar[d] \\
Y \ar[rr] & & \spec\K\,.}
\end{equation}

(2) By  definition
\begin{equation*}
(\s{F}_2\lboxtimes\s{F}_1)\FMcomp\c W=
\R{\pr_{1,3}}_*\left(\pr_{2,3}^*(\pr_2^*\s{F}_2\lotimes\pr_1^*\s{F}_1)
\lotimes\pr_{1,2}^*\c W\right).
\end{equation*}
Since $\pr_2\comp \pr_{2,3}=\pr_2\comp\pr_{1,3}$ and $\pr_1\comp
\pr_{2,3}=\pr_2\comp\pr_{1,2}$, the last expression can be rewritten as
\begin{equation*}
\R{\pr_{1,3}}_*\left(\pr_{1,3}^*\pr_2^*{\s{F}_2}\lotimes
\pr_{1,2}^*(\pr_2^*{\s{F}_1} \lotimes\c W)\right).
\end{equation*}
Using the projection formula and the ``flat base change'' theorem for
the Cartesian square
\begin{equation}\label{e33}
\xymatrix{
Y\Times Y\Times Y \ar[rr]^{\pr_{1,2}} \ar[d]_{\pr_{1,3}} & &
Y\Times Y \ar[d]^{\pr_1 }  \\
Y\Times Y \ar[rr] ^{ \pr_1} & & Y}
\end{equation}
we arrive at
\begin{equation*}
\pr_2^*\s{F}_2\lotimes
\R{\pr_{1,3}}_*\pr_{1,2}^*(\pr_2^*\s{F}_1\lotimes\c W) \iso
\pr_2^*\s{F}_2\lotimes\pr_1^*\R{\pr_1}_*(\pr_2^*\s{F}_1\lotimes\c W)=
{\s{F}_2}\lboxtimes\FM[2]{\c W}({\s{F}_1}).
\end{equation*}
This proves (2).

The proof of (3) is completely similar to the proof of (2), while (4)
follows immediately from the first two statements.
\end{proof}

Let $\diag\colon Y\to Y\Times Y$ be the diagonal morphism.
We will write $\so_{\Delta_Y}$ (or simply $\so_{\Delta}$) for
$\diag_*\so_Y$. It is well known that
$\FM[i]{\so_{\Delta}}\iso\id_{\D^b(Y)}$ for $i=1,2$, and more
generally, if $\s L$ is a line bundle on $Y$, $\FM[i]{\diag_*\s L}$ is
isomorphic to the autoequivalence of $\D^b(Y)$ defined by $\s
F\mapsto\s F\otimes\s L$ (see, e.g., \cite[Ex. 5.4]{Huybrechts}).

\section{The resolution of the diagonal for weighted projective spaces}
\label{s:WB}

Let $\P:=\P^n({\bf w})$ be a weighted projective space (regarded as a
stack)\footnote{We will often refer to \cite{Alberto}, where weighted
projective spaces are regarded as graded schemes. The equivalence with
the point of view of stacks is explained in \cite[\S 1.6]{Alberto}}
with weight vector ${\bf w}= (\w_0,\w_1,\ldots,\w_n)$,
and $\w_i>0$ for all $i=0,\ldots,n$. We introduce some notations that will be
used throughout the paper. For a subset $I\subseteq\{0,1,\ldots,n\}$ the symbol
$\card{I}$ denotes the cardinality of $I$. Similarly we introduce the following
sums of weights
\begin{equation*}
\w_I = \sum_{i\in I}\w_i\,;\qquad \sww =
\w_{\{0,1,\ldots,n\}}=\sum_{i=0}^n\w_i\,.
\end{equation*}

Let $\p=\K[x_0,x_1,\ldots,x_n]$ be the graded polynomial ring where the
generators have $\deg(x_i)=\w_i$. First recall the Koszul complex $\sko$ on
$\P=\P^n({\bf w})$ associated to the regular sequence $(x_0,x_1,\ldots,x_n)$,
whose $j$th term is given by
\[
\sko^j=\bigoplus_{\card{I}=-j}\so(-\w_I)\,,
\]
where for a given $i\in\Z$ we abbreviated $\so_{\P}(i)$ by $\so(i)$. The
summation is over all subsets of $\{0,1,\ldots,n\}$ of cardinality $-j$.
Obviously, $\sko^j$ is non-zero only for $-n-1\le j\le0$. The components of the
$j$th differential of $\sko$
\[
\xymatrix@1{\ar[rr]^-{\diff{\sko}^j} \sko^j = \bigoplus_{\card{I}=-j}\so(-\w_I)
\, &&\, \sko^{j+1} = \bigoplus_{\card{I'}=-j-1}\so(-\w_{I'})}
\]
are given by
\begin{equation*}
\left(\diff{\sko}^j\right)^I_{I'}=\begin{cases}
(-1)^{N_I^i}\,x_i\in \p_{\w_i} & \text{if $I=I'\cup\{i\}$} \\
0 & \text{otherwise}
\end{cases}
\end{equation*}
where the integer $N_I^i$ is defined as the cardinality
\begin{equation*}
N_I^i = \card{\{j\in I\colon j<i\}}.
\end{equation*}
The notation $\p_a$, for $a\in\Z$, refers to the degree $a$ subspace of
$\p=\K[x_0,x_1,\ldots,x_n]$; $\p_a$ will be viewed as a $\K$-vector space.
Since we chose the weights of the $x_i$'s to be positive, $\p_a$ is
the zero vector space whenever $a<0$.

Following \cite[Definition 2.5.2]{Alberto} for $-\sww<l\le0$ we introduce the
subcomplex $\bo{l}$ of the twisted Koszul complex $\sko(-l)$:
\begin{equation*}
\bo{l}^j = \bigoplus_{\card{I}=-j,\w_I\leq -l}\so(-l-\w_I)\quad\subseteq\quad
\bigoplus_{\card{I}=-j}\so(-l-\w_I) = \sko^j(-l)\,.
\end{equation*}
Note that $\bo{0}$ is the complex $\so$ concentrated in degree 0.

The $\bo{l}$'s have a natural interpretation: the left dual of the full
and strong exceptional sequence $(\so,\so(1),\ldots,\so(\sww-1))$ in $\D^b(\P)$
is the full exceptional sequence
$(\bo{1-\sww}[1-\sww],\ldots,\bo{-1}[-1],\bo{0})$
\cite[Proposition 2.5.11]{Alberto}.
We will use the complexes $\bo{l}$ to give a generalization of Beilinson
resolution of the diagonal for the stack $\P$.

For every $-\sww<k\le0$ we define the complex $\pdr{k}\in\CC^b(\P\Times\P)$
inductively. The starting point is
\begin{equation*}
\pdr{0}:=\so_{\P}\boxtimes{\bo{0}}\iso\so_{\P\times \P}.
\end{equation*}
For $-\sww<k<0$ we set
\begin{equation}\label{mcdiag}
\pdr{k}:= \mc{\mpdr{k}\colon\so(k)\boxtimes\bo{k}[-1]\too\pdr{k+1}},
\end{equation}
where $\mpdr{k}$ is a natural map that will be defined below.
Observe that $\pdr{k}$ is defined in terms of $\pdr{k+1}$ since $k$ is
increasingly more negative. This convention ties well with the fact that all of
our complexes will be non-trivial {\em only in negative} degree.

Given the definition of the mapping cone in \eqref{vv1}, for a fixed
$k$, the components of the
complex $\pdr{k}$ are immediate: for $j\in\Z$ they are
\begin{equation}\label{e2}
\pdr{k}^j=\bigoplus_{k\le l\le0}\so(l)\boxtimes\bo{l}^j=
\bigoplus_{\substack{k\le l\le0\\ \card{I}=-j,\w_I\le-l}}\so(l,-l-\w_I).
\end{equation}
Here we use the shorthand notation
\[
\so(i,j) = \so_{\P\times \P}(i,j):=\so_{\P}(i)\boxtimes\so_{\P}(j),\qquad
\txt{for $i,j\in\Z$ }\,.
\]

As a result, the $j$th component of $\mpdr{k}$ in \eqref{mcdiag} has to be a
map between
\begin{equation}\label{e7}
\mpdr{k}^j\colon\bigoplus_{\card{I}=1-j,\w_I\le-k}\so(k,-k-\w_I)
\too\bigoplus_{\substack{k<l\le0\\ \card{I}=-j,\w_I\le-l}}\so(l,-l-\w_I).
\end{equation}
Now each component of $\mpdr{k}^j$
\begin{equation*}
\left(\mpdr{k}^j\right)^I_{l,I'}\in
\Hom_{\P\times \P}(\so(k,-k-\w_I),\so(l,-l-\w_{I'}))\iso
\p_{l-k}\otimes_{\K}\p_{k-l+\w_I-\w_{I'}}
\end{equation*}
for $k<l\le0$, $\card{I}=1-j$, $\card{I'}=-j$, $\w_I\le-k$ and
$\w_{I'}\le-l$, is defined to be
\begin{equation}\label{e6}
\left(\mpdr{k}^j\right)^I_{l,I'}:=\begin{cases}
-(-1)^{N_I^i}\,x_i\otimes 1\in\p_{\w_i}\otimes_{\K}\p_0 &
\text{if $l=k+\w_i$, $I=I'\cup\{i\}$} \\
0 & \text{otherwise}.
\end{cases}
\end{equation}
Note that the map $x_i\otimes 1$ imposes the conditions $l=k+\w_i$ and
$-k-\w_I=-l-\w_{I'}$, which lead to $\w_I=\w_{I'}+\w_i$, and is
automatically satisfied by the condition $I=I'\cup\{i\}$.

For the inductive definition in \eqref{mcdiag} to make sense we need to
make sure that for any $-\sww<k<0$
\begin{enumerate}
 \item[(i.)] $\pdr{k+1}$ is a complex,
 \item[(ii.)] $\mpdr{k}$ is  a morphism of complexes.
\end{enumerate}
We can prove these simultaneously by induction: $\pdr{0}=
\so_{\P}\boxtimes\bo{0}$ is clearly a complex; by assuming that $\pdr{k+1}$ is
a complex and that $\mpdr{k}$ is  a map of complexes, the mapping cone
construction in \eqref{mcdiag} guarantees that $\pdr{k}$ is a complex, i.e.,
$\diff{\pdr{k}}$ is a differential. Therefore, the key point is to prove (ii.),
and then (i.) follows automatically.

In the light of \eqref{e2} and \eqref{e6} and the definition \eqref{mcdiag},
it is immediate to write down
explicitly the candidate differentials of $\pdr{k}$ (so far these are
only maps, since we have not yet proven that $\diff{\pdr{k}}^{\;2}=0$). Each
component
\[
\left(\diff{\pdr{k}}^j\right)^{l,I}_{l',I'}\in
\Hom_{\P\times \P}(\so(l,-l-\w_I),\so(l',-l'-\w_{I'}))\iso
\p_{l'-l}\otimes_{\K}\p_{l-l'+\w_I-\w_{I'}}
\]
of $\diff{\pdr{k}}^j\colon\pdr{k}^j\to\pdr{k}^{j+1}$, for $k\le
l,l'\le0$, $\card{I}=-j$, $\card{I'}=-j-1$, $\w_I\le-l$ and
$\w_{I'}\le-l'$, is given by
\begin{equation}\label{e1}
\left(\diff{\pdr{k}}^j\right)^{l,I}_{l',I'}=
\begin{cases}
\phantom{-}(-1)^{N_I^i}\,1\otimes x_i\in\p_0\otimes_{\K}\p_{\w_i} &
\text{if $l'=l$, $I=I'\cup\{i\}$} \\
-(-1)^{N_I^i}\,x_i\otimes 1\in\p_{\w_i}\otimes_{\K}\p_0 &
\text{if $l'=l+\w_i$, $I=I'\cup\{i\}$} \\
0 & \text{otherwise}.
\end{cases}
\end{equation}
This is straightforward to check using the definition of mapping cone
and assuming inductively that it holds for
$k+1$. Thus everything will be well-defined once we  prove the
following:

\begin{lemma}
$\mpdr{k}\colon\so(k)\boxtimes\bo{k}[-1]\too\pdr{k+1}$ as defined above is a
map of complexes.
\end{lemma}

\begin{proof}
By the inductive hypothesis discussed above,
we can assume that $\pdr{k+1}$ is a complex,
and that $\diff{\pdr{k+1}}$ is given by \eqref{e1}.
Glancing at the definition \eqref{mcdiag}, we need to show that
\[
\mpdr{k}^{j+1}\comp\diff{\so(k)\boxtimes\bo{k}[-1]}^j=
\diff{\pdr{k+1}}^j\comp\mpdr{k}^j.
\]
Using the fact that $\diff{\so(k)\boxtimes{\bo{k}}[-1]}^j=
-\id_{\so(k)}\boxtimes\,\diff{\bo{k}}^{j-1}$, this is equivalent
to the following square being commutative
\begin{equation}\label{e8}
\xymatrix{\so(k)\boxtimes\bo{k}^{j-1}
\ar[rrr]^{-\id_{\so(k)}\boxtimes\,\diff{\bo{k}}^{j-1}}
\ar[d]^{\mpdr{k}^j} & & & \so(k)\boxtimes\bo{k}^j
\ar[d]^{\mpdr{k}^{j+1}} \\
\pdr{k+1}^j \ar[rrr]^{\diff{\pdr{k+1}}^j} & & & \pdr{k+1}^{j+1}\,.}
\end{equation}
Let us restrict to one of the direct summands of
$\so(k)\boxtimes\bo{k}^{j-1}$, say
$\so(k,-k-\w_I)$. The square \eqref{e8} involves two mappings: the
``horizontal then vertical'' map is
\begin{equation*}
\xymatrix{\so(k,-k-\w_I)
\ar[rrrrr]^-{\bigoplus_{\substack{i\in I\\ i'\in I_i}}\,
(-1)^{N^i_I+N^{i'}_{I_i}}(x_{i'}\otimes x_i)} & & & & &
\bigoplus_{\substack{i,i'\in I\\ i\ne i'}}\,\so(k+\w_{i'},\w_i-k-\w_I)}
\end{equation*}
where $I_i=I-\{i\}$; while the ``vertical then horizontal'' map is
\begin{equation*}
\xymatrix{\so(k,-k-\w_I)
\ar[rrrrrr]^-{\bigoplus_{\substack{i\in I\\ i'\in I_i}}\,
-(-1)^{N^i_I+N^{i'}_{I_i}}\txt{$(x_i\otimes x_{i'})\oplus$\\
$(-x_ix_{i'}\otimes1)$}} & & & & & &
\bigoplus_{\substack{i,i'\in I\\ i\ne i'}}\,
\txt{$\so(k+\w_i,\w_{i'}-k-\w_I)\oplus$\\
$\so(k+\w_i+\w_{i'},-k-\w_I)$}\,.}
\end{equation*}
Let us focus on the second map. Writing $\bigoplus_{i\in I,i'\in I_i}$
artificially singled out one element of the set $\{i,i'\}$, since
the summation is over pairs $i\ne i'$. Changing variables
$i\leftrightharpoons i'$ and observing that
\[
(-1)^{N^i_I+N^{i'}_{I_i}}=-(-1)^{N^{i'}_I+N^i_{I_{i'}}}
\]
we see immediately that the second component,
$(-x_ix_{i'}\otimes1)$, is
in fact the zero map; whereas after the exchange $i\leftrightharpoons
i'$ the first component is identical to the ``horizontal then
vertical'' map. This proves that the square indeed commutes.
\end{proof}

The main benefit of these definitions is the following generalization
to weighted projective spaces of Beilinson's resolution of the
diagonal for $\P^n$:

\begin{prop}\label{diagres}
There is a natural morphism of complexes
$\ndr\colon\pdr{1-\sww}\to\so_{\Delta}$,
which descends to an isomorphism in $\D^b(\P\Times\P)$.
\end{prop}

\begin{proof}
Let us set for brevity $\dr:=\pdr{1-\sww}$. In order to construct
the natural morphism of complexes $\ndr\colon\dr\to\so_{\Delta}$,
start with the adjunction $\diag^*\dashv\diag_*$ and observe that
\[
\Hom_{\P\times \P}(\so(l,-l),\so_{\Delta}=\diag_*\so_{\P})\,\iso\,
\Hom_{\P}(\diag^*\so(l,-l),\so_{\P})\iso\Hom_{\P}(\so_{\P},\so_{\P})
\]
for any integer $l$. Let $f_l\colon\so(l,-l)\to\so_{\Delta}$ be the
morphism corresponding to the identity under this isomorphism.
Equivalently, $f_l$ is induced by ``multiplication''. This follows
from the fact that
\[
\so(l,-l)=\pr_2^*{\so(l)}\otimes\pr_1^*{\so(-l)}
\iso\pr_2^*{\so(l)}\otimes\pr_1^*{\so(l)\dual}
\]
and we can use the natural pairing between $\so(l)$ and $\so(l)\dual$ to map
$\so(l,-l)$ into $\so_{\Delta}$.

Since $\dr^0=\bigoplus_{-\sww<l\le0}\so(l,-l)$, we define $\ndr^0:=
\bigoplus_{-\sww<l\le0}f_l$; while of course for $i\neq 0$ we let
$\ndr^i=0$. To check that $\ndr\colon\dr\to\so_{\Delta}$ is a morphism
of complexes, we only need to show that
$\ndr^0\comp\diff{\dr}^{-1}=0$. The
relevant diagram is
\begin{equation*}
\xymatrix{\dr \ar[d]^{\ndr\quad:} & &
\displaystyle{\bigoplus_{-\sww<l\le0,\w_i\le-l}}\so(l,-l-\w_i)
\ar[rr]^-{\diff{\dr}^{-1}} \ar[d]^{\ndr^{-1}=0} & &
\displaystyle{\bigoplus_{-\sww<l\le0}}\so(l,-l) \ar[d]^{\ndr^0} \\
\so_{\Delta} & & 0 \ar[rr] & & \so_{\Delta}.}
\end{equation*}
Pick a component in $\dr^{-1}$, say $\so(l,-l-\w_i)$. There are two
non-trivial components of $\diff{\dr}^{-1}$ emanating from it:
$\left(\diff{\dr}^{-1}\right)^{l,\{i\}}_{l,\emptyset}$ and
$\left(\diff{\dr}^{-1}\right)^{l,\{i\}}_{l+\w_i,\emptyset}$; and it is
clear from \eqref{e1} that
\[
f_l\,\comp\left(\diff{\dr}^{-1}\right)^{l,\{i\}}_{l,\emptyset}+f_{l+\w_i}
\,\comp\left(\diff{\dr}^{-1}\right)^{l,\{i\}}_{l+\w_i,\emptyset}=0.
\]
This proves that $\ndr\colon\dr\to\so_{\Delta}$ is indeed a morphism of
complexes.

To show that $\ndr$ is an isomorphism in $\D^b(\P\Times\P)$ it suffices to
prove that $\FM[1]{\ndr}\colon\FM[1]{\dr}\to\FM[1]{\so_{\Delta}}$ is an
isomorphism of functors. Indeed, assuming that $\FM[1]{\ndr}$ is an
isomorphism, we immediately deduce the isomorphism of functors
$\FM[1]{\coneI{\ndr}}\iso0$, which implies that $\cone{\ndr}\iso0$ (hence
$\ndr$ is an isomorphism) in $\D^b(\P\Times\P)$ by Lemma~\ref{ker0}.

Now recall that $(\so(1-\sww),\ldots,\so(-1),\so)$ is a full and
strong exceptional sequence (\cite[Remark 2.2.6]{Alberto}, \cite[Theorem
2.12]{Auroux}), and hence, in particular, it generates $\D^b(\P)$ as
triangulated category. 
Therefore, in order to prove that $\FM[1]{\ndr}$ is an
isomorphism it is enough to show that
\[
\FM[1]{\ndr}(\so(k))\colon\FM[1]{\dr}(\so(k))\too\FM[1]{\so_{\Delta}}(\so(k))
\]
is an isomorphism for $-\sww<k\le0$.

It is immediate that $\FM[1]{\so_{\Delta}}(\so(k))\iso\so(k)$. Next we
compute $\FM[1]{\dr}(\so(k))$ in two different ways: first using the
recursion \eqref{mcdiag}, and then using \eqref{e2}. For the first
computation
notice that by part (1) of Lemma~\ref{kercomp}
\[\FM[1]{\so(l)\boxtimes\bo{l}}(\so(k))\iso
\so(l)\otimes_{\K}\R\Gamma(\P,\bo{l}\otimes\so(k)).\]
On the other hand, for $-\sww<k,l\le0$, $\dim
H^i(\P,\bo{l}(k))=\delta_{k,l}\delta_{i,0}$, as shown in the proof of
\cite[Theorem 2.5.8]{Alberto}; therefore
\[\FM[1]{\so(l)\boxtimes\bo{l}}(\so(k))\iso
\begin{cases}
\so(k) & \text{if $k=l$} \\
0 & \text{if $k\ne l$}.
\end{cases}\]
It is also immediate that
\[\FM[1]{\mcI{f\colon\s{A}\to\s{B}}}(\s{F})\iso
\cone{\FM[1]{f}(\s{F})\colon\FM[1]{\s{A}}(\s{F})\too\FM[1]{\s{B}}(\s{F})}.\]
Using these two facts and the defining equation \eqref{mcdiag},
it is easy to deduce that
\begin{equation}\label{e5}
\FM[1]{\dr}(\so(k))\iso\so(k).
\end{equation}

At this point it is clear that $\FM[1]{\ndr}(\so(k))$ can be
identified with a map $\so(k)\to\so(k)$, but it is not evident which map
it is. To settle this question we can also compute
$\FM[1]{\dr}(\so(k))$ directly from \eqref{e2} and \eqref{e1}. Observe
that by the projection formula
\begin{equation*}
\R{\pr_2}_*\so(i,j)=
\R{\pr_2}_*(\pr_2^*\so(i)\otimes\pr_1^*\so(j))\iso
\so(i)\otimes\R{\pr_2}_*\pr_1^*\so(j).
\end{equation*}
On the other hand, for any $j>-\sww$, from \eqref{e22}
we know that $\R{\pr_2}_*\pr_1^*\so(j)=\so\otimes_{\K}\p_j;$
and as a result
\begin{equation}\label{f2}
\R{\pr_2}_*\so(i,j)\iso{\pr_2}_*\so(i,j)\iso\so(i)\otimes_{\K}\p_j.
\end{equation}
Therefore, for any $-\sww<k\le0$, every term of
$\dr\otimes\pr_2^*\so(k)$ is ${\pr_1}_*$-acyclic, and we have a
natural isomorphism $\FM[1]{\dr}(\so(k))\iso
{\pr_1}_*(\dr\otimes\pr_2^*\so(k))\iso\drim{k}$, where
\[
\drim{k}^j=\bigoplus_{-\sww<l\le0,\card{I}=-j,\w_I\le-l}
\so(l)\otimes_{\K}\p_{k-l-\w_I}
\]
and $\diff{\drim{k}}^j$ is induced by $\diff{\dr}^j$,
and in fact is given by exactly the same
formal expression as \eqref{e1}, although the two act
on different complexes. Also note that most
terms in $\drim{k}^j$ are zero, except for those that satisfy
the condition $k-l-\w_I\geq 0$.

In this presentation $\FM[1]{\ndr}(\so(k))$ can be identified with the
natural morphism $\mdrim{k}\colon\drim{k}\to\so(k)$, where
the only non-zero component of $\mdrim{k}$ is
\[
\mdrim{k}^0: \bigoplus_{-\sww<l\le0}\so(l)\otimes_{\K}\p_{k-l}\too\so(k),
\]
and this is such that each $\mdrim{k}^0\rest{\so(l)\otimes_{\K}\p_{k-l}}$ is
given by the multiplication map. This follows readily from the definition of
$\ndr\colon\dr\to\so_{\Delta}$ and the isomorphisms that led to
$\FM[1]{\dr}(\so(k))\iso\drim{k}$.

Clearly $\mdrim{k}^0$ is surjective, hence
\[
H^0(\mdrim{k}):\; H^0(\drim{k})=\drim{k}^0\, /\im\diff{\drim{k}}^{-1}\too\so(k)
\]
is also surjective. We have already seen in \eqref{e5} that
$H^0(\drim{k})\iso\so(k)$ and $H^i(\drim{k})=0$ for $i\ne0$, and thus
$H^0(\mdrim{k})$ is an isomorphism; and consequently
$\mdrim{k}$ is a quasi-isomorphism.\footnote{In fact $\mdrim{k}$ is
an isomorphism in $\KK^b(\P)$, with inverse given by the natural morphism
$\so(k)\mono\drim{k}^0$, but this statement is of marginal interest to us.}
\end{proof}

\begin{remark}
{\em 
If $\w_0=\cdots=\w_n=1$, i.e., $\P$ is the ordinary projective space $\P^n$,
then the complex $\dr=\pdr{-n}$ defined in \eqref{mcdiag}
does not coincide with the well-known resolution of the diagonal
first considered by Beilinson \cite{Bei:res} (for weighted
projective planes, i.e., $n=2$, the resolution coincides with the one
considered in the context of quivers by King \cite{King}). 
Beilinson's resolution
$\bdr\in\CC^b(\P^n\Times\P^n)$ is defined by $\bdr^j:=
\so(j)\boxtimes\sd^{-j}(-j)$ (again, $\bdr^j\ne0$ only for $-n\le
j\le0$) with differential given (for $-n\le j<0$) by the natural morphism
\begin{multline*}
\diff{\bdr}^j\in\Hom_{\P^n\Timesa\P^n}
(\so(j)\boxtimes\sd^{-j}(-j),\so(j+1)\boxtimes\sd^{-j-1}(-j-1))\\
\iso\Hom_{\P^n}(\so(j),\so(j+1))\otimes_{\K}
\Hom_{\P^n}(\sd^{-j}(-j),\sd^{-j-1}(-j-1))\iso
\p_1\otimes_{\K}\p_1\dual\iso\Hom_{\K}(\p_1,\p_1)
\end{multline*}
corresponding to $\id$. Observe that if we define complexes
$\bpdr{k}$ by
\[\bpdr{k}^j:=\begin{cases}
\bdr^j & \text{if $j\ge k$} \\
0 & \text{if $j<k$}
\end{cases}\]
with $\diff{\bpdr{k}}^j:=\diff{\bdr}^j$ for $j\ge k$ (hence
$\bdr=\bpdr{-n}$), then $\bpdr{0}=\so_{\P^n\Times\P^n}=\pdr{0}$, while for
$-n\le k<0$
\[\bpdr{k}=\mc{\so(k)\boxtimes\sd^{-k}(-k)[-k-1]\too\bpdr{k+1}}\]
(where the morphism is induced by $\diff{\bdr}^k$), in analogy with
\eqref{mcdiag}. Actually, using
the fact that $\bo{k}\iso\sd^{-k}(-k)[-k]$ in $\D^b(\P^n)$ (see
\cite[Remark 2.5.9]{Alberto}), it can be easily proved by descending
induction on $k$ that $\bpdr{k}\iso\pdr{k}$ in $\D^b(\P^n\Times\P^n)$
(but not in $\KK^b(\P^n\Times\P^n)$). Clearly the complex $\bdr$ is
simpler than $\dr$, but it does not seem possible to extend it to the
weighted case, the problem coming from the fact that the
complexes $\bo{k}$ are not quasi-isomorphic to shifts of ordinary
sheaves in general.
}
\end{remark}

\section{Proof of the theorem}    \label{s:proof}

We start out by expressing the functors $\ms{K}$ and $\ms{L}$ appearing in
Theorem~\ref{t1}, and defined in \eqref{k2} and \eqref{k3}, as Fourier-Mukai
functors. Lemma 3.2 of \cite{ST:braid} shows that the kernel of $\ms{K}$ is $\c
K:=\cone{\so_{X\Timesa X}\overset{\diag\mrs}{\too}\diag_*\so_X}$, where
$\diag\mrs$ is the natural map associated to $\diag\colon X\mono X\Times
X$ (clearly $\c K\iso\c I_{\Delta}[1]$, where $\c I_{\Delta}$ is the
ideal sheaf of the diagonal in $X\Times X$, but this observation is of
no use in our context).
Similarly, the kernel of $\ms{L}$ is $\c L:=\diag_*\so_X(1)$.
Therefore $\ms{G}=\ms{L} \,\comp \ms{K} \iso\FM[1]{\s{L}\FMcomp\s{K}}$. Using
(3) of Lemma~\ref{kercomp} it is easy to show that
\[
 \s{L}\FMcomp\s{K}=\diag_*\so_X(1) \FMcomp
 \cone{\so_{X\Timesa X}\overset{\diag\mrs}{\too}\diag_*\so_X}\iso
\cone{\so_{X\Timesa X}(1,0)\overset{g}{\too}\diag_*\so_X(1)},
\]
where $g\colon\so_{X\Timesa X}(1,0)\to\diag_*\so_X(1)$ is the natural morphism.
Using this fact, Theorem~\ref{t1} is a consequence of the following
proposition:

\begin{prop}\label{mainprop}
Let $X$ be a smooth anti-canonical stacky hypersurface in $\P$. For the natural
morphism $g\colon\so_{X\Timesa X}(1,0)\to \diag_*\so_X(1)$ let $\G$ denote its
mapping cone; i.e, $\G:=\mc{g}$. Then $(\G)^{\FMcomp\sww}\iso\so_{\Delta_X}[2]$
in $\D^b(X\Times X)$.
\end{prop}

The rest of the section is dedicated to the proof of the proposition.

Let $\iota\colon X\mono\P$ denote the inclusion. For a complex $\s{E}
\in\KK^b(\P)$ we will write $\cprest{\s{E}}$ for the restriction
$\iota^*{\s{E}}\in\KK^b(X)$ and similarly for morphisms in $\KK^b(\P)$.
The same notation will be used for the inclusion
$\iota\Times\iota\colon X\Times X\mono\P\Times\P$, but it will be
clear from the context which one is meant. Note that if each $\s{E}^j$
is locally free, then $\cprest{\s{E}}\iso\Ll\iota^*{\s{E}}$ in $\D^b(X)$.

For $0<m\le\sww$ we define the complex
\begin{equation}\label{f1}
\Gi{m}:=\mc{\pdrX{1-m}\overset{\zeta_m}{\too}\diag_*\so_X}
\end{equation}
where each component of the morphism
\begin{equation}\label{f17}
\zeta_m^0\colon\pdrX{1-m}^0=
\bigoplus_{-m<l\le0}\so_{X\Timesa X}(l,-l)\to\diag_*\so_X
\end{equation}
is the
natural map induced by multiplication (the argument in the
proof of Prop.~\ref{diagres}  showing that
$\ndr$ is a morphism of complexes also shows that $\zeta_m$ is a
morphism of complexes).

We claim that it suffices to prove that
\begin{claim}\label{e3}
$(\G)^{\FMcomp m}\iso{\Gi{m}}(m,0)$ for all $\,0<m\le\sww$.
\end{claim}

Indeed, assuming this, and taking into account
that $(\diag_*\so_X)(m,0)\iso\diag_*\so_X(m)$, the proposition follows
from the following:

\begin{lemma}\label{l0}
$\Gi{\sww}\iso\diag_*\so_X(-\sww)[2]$ in $\D^b(X\Times X)$.
\end{lemma}

\begin{proof}
Since $\iota\Times\iota$ is a closed immersion, it is sufficient to prove that
\begin{equation}\label{e10}
(\iota\Times\iota)_*\Gi{\sww}\iso(\iota\Times\iota)_*\diag_*\so_X(-\sww)[2]
\end{equation}
in $\D^b(\P\Times\P)$. By the projection formula and using the fact
that $\diag\comp \iota =(\iota\Times\iota) \comp\diag$ (the first
$\diag$ is the diagonal of $\P$, while the second is the diagonal of $X$),
\[
(\iota\Times\iota)_*\Gi{\sww}\iso
\cone{(\iota\Times\iota)_*\pdrX{1-\sww}\too(\iota\Times\iota)_*\diag_*\so_X}
\iso\cone{\pdr{1-\sww}\otimes(\iota\Times\iota)_*\so_{X\Timesa X}\too
\diag_*\iota_*\so_X}.
\]
It is useful to replace
$(\iota\Times\iota)_*\so_{X\Timesa X}$ with a locally free resolution. To this
purpose, let us start with the short exact sequence defining $X$
\begin{equation}\label{eX}
\xymatrix@1{0 \ar[r] & \so_{\P}(-\sww) \ar[r]^-{s} & \so_{\P}
\ar[r]^-{\iota\mrs} & \iota_*\so_X \ar[r] & 0}.
\end{equation}

Since we want a resolution for $(\iota\Times\iota)_*\so_{X\Timesa X}$
we consider the complex $\s{S}$:
\[
\s{S}:\qquad
0 \to\so_{\P\times \P}(-\sww,-\sww)\mor{\begin{pmatrix}
-1\otimes s \\
s\otimes 1
\end{pmatrix}}
\begin{matrix}
\so_{\P\times \P}(-\sww,0) \\
\oplus \\
\so_{\P\times \P}(0,-\sww)
\end{matrix}\mor{\begin{pmatrix}
s\otimes1 & 1\otimes s
\end{pmatrix}}\so_{\P\times \P}\to0
\]
where $\so_{\P\times \P}$ is sitting in degree 0.

As $X\Times X$ is a codimension two complete intersection in
$\P\Times\P$, defined precisely by $s\otimes1$ and $1\otimes s$, it is
clear that the canonical map
$(\iota\Times\iota)\mrs\colon\s{S}^0=\so_{\P\times \P}\to
(\iota\Times\iota)_*\so_{X\Timesa X}$
induces a quasi-isomorphism 
between $\s{S}$ and $(\iota\Times\iota)_*\so_{X\Timesa X}$. Thus
\[
(\iota\Times\iota)_*\Gi{\sww}\iso
\cone{f\colon\pdr{1-\sww}\otimes\s{S}\too\diag_*\iota_*\so_X}.
\]
In the light of \eqref{f17} each
component of $f^0\colon(\pdr{1-\sww}\otimes\s{S})^0\iso
\bigoplus_{-\sww<l\le0}\so_{\P\times \P}(l,-l)\to\diag_*\iota_*\so_X$ is
the natural one (corresponding to $\iota\mrs$ under adjunction).

Observing that $\diag^*\s{S}$ can be identified with the complex
\[0 \to\so_{\P}(-2\sww)\mor{\begin{pmatrix}
-s \\
s
\end{pmatrix}}
\begin{matrix}
\so_{\P}(-\sww) \\
\oplus \\
\so_{\P}(-\sww)
\end{matrix}\mor{\begin{pmatrix}
s & s
\end{pmatrix}}\so_{\P}\to0\] and denoting by
$h\colon\diag^*\s{S}\to\iota_*\so_X$ the natural morphism defined by
$h^0=\iota\mrs\colon(\diag^*\s{S})^0\iso\so_{\P}\to\iota_*\so_X$, it is also
straightforward to check that the diagram
\[
\xymatrix@C=15mm{
\pdr{1-\sww}\otimes\s{S} \ar[r]^-f \ar[d]_{\ndr\otimes\id} &
\diag_*\iota_*\so_X \\
(\diag_*\so_{\P})\otimes\s{S} \ar[r]_-{\iso}^-{\varpi} & \diag_*\diag^*\s{S}
\ar[u]^{\diag_*h}}
\]
commutes ($\ndr\colon\pdr{1-\sww}\to\so_{\Delta}$ is the morphism from
Prop.~\ref{diagres}). Using this and Prop.~\ref{diagres} we have that
\[
(\iota\Times\iota)_*\Gi{\sww}\iso\cone{f}\iso
\cone{(\diag_*h)\comp\varpi\comp(\ndr\otimes\id)}
\iso\cone{\diag_*h}\iso\diag_*\cone{h}.
\]
On the other hand $\diag^*\s{S}$ is isomorphic in $\CC^b(\P)$
to the complex
\[
\diag^*\s{S}:\qquad
0 \to\so_{\P}(-2\sww)\mor{\begin{pmatrix}
0 \\
s
\end{pmatrix}}
\begin{matrix}
\so_{\P}(-\sww) \\
\oplus \\
\so_{\P}(-\sww)
\end{matrix}\mor{\begin{pmatrix}
s & 0
\end{pmatrix}}\so_{\P}\to0,\]
which clearly splits as $\mc{\so_{\P}(-\sww)\overset{s}{\too}\so_{\P}}\oplus
\mc{\so_{\P}(-2\sww)\overset{s}{\too}\so_{\P}(-\sww)}[1]$.
Since we defined $h^0$ to be $\iota\mrs\colon \so_{\P}\to\iota_*\so_X$
we have that
\[
 \cone{h}\iso \left(\so_{\P}(-\sww)\overset{s}{\too}
 \so_{\P}\overset{\iota\mrs}{\too}
\iota_*\so_X\right)\oplus
\mc{\so_{\P}(-2\sww)\overset{s}{\too}\so_{\P}(-\sww)}[2].
\]
Taking into account \eqref{eX} it is clear that
$\cone{h}\iso\iota_*\so_X(-\sww)[2]$ in $\D^b(\P)$. Therefore
\[(\iota\Times\iota)_*\Gi{\sww}\iso\diag_*\cone{h}\iso
\diag_*\iota_*\so_X(-\sww)[2]\iso(\iota\Times\iota)_*\diag_*\so_X(-\sww)[2],\]
which proves \eqref{e10}, and concludes the proof of the lemma.
\end{proof}

Returning to the proof of the proposition,
we will prove Claim~\ref{e3} by induction on $m$.
For $m=1$ the isomorphism $\G
\iso\Gi{1}(1,0)$ immediately follows from the definitions of
$\G$ and $\Gi{1}$. Therefore we assume
Claim~\ref{e3} to hold for some $0<m<\sww$, and prove it for $m+1$.

Noticing that $\G\iso\diag_*\so_X(1)\FMcomp\Gi{1}$, what we need to
show is equivalent to
\begin{equation}\label{e0}
\Gi{1}\FMcomp\Gi{m}(m,0)\iso\Gi{m+1}(m,0).
\end{equation}
In order to prove this, we start with the following lemma:

\begin{lemma}\label{l1}
For any $\s{E}\in\D^b(X\Times X)$ we have a natural isomorphism
\[
\Gi{1}\FMcomp\s{E}\iso
\cone{\so_X\boxtimes \R{\pr_1}_*{\s{E}} \too\s{E}}\iso
\cone{\pr_1^*\R{\pr_1}_*{\s{E}}\too\s{E}},
\]
where the morphism in the right hand side is the natural one
(corresponding, under adjunction, to $\id_{\R{\pr_1}_*\s{E}}$).
\end{lemma}

\begin{proof}
If $f\colon Y\to Z$ is a morphism of stacks,
$\lradj{f}\colon\Ll f^*\comp\R f_*\to\id_{\D^b(Y)}$ and
$\rladj{f}\colon\id_{\D^b(Z)}\to\R f_*\comp\Ll f^*$ will denote the
adjunction morphisms; then we need to prove that
$\Gi{1}\FMcomp\s{E}\iso\cone{\lradj{\pr_1}(\s{E})}$. Since
$(-)\FMcomp\s{E}$ is an exact functor we have
$\Gi{1}\FMcomp\s{E}\iso\cone{\diag\mrs\FMcomp\id_{\s{E}}\colon
\so_{X\Timesa X}\FMcomp\s{E}\to\diag_*\so_X\FMcomp\s{E}}$.
Applying the ``flat base change'' theorem to the Cartesian square
\[\xymatrix{
X\Times X \ar[rr]^{\tilde{\diag}} \ar[d]_{\pr_2} & &
X\Times X\Times X \ar[d]^{\pr_{2,3}}  \\
X \ar[rr]^{\diag} & & X\Times X
}\]
and using the projection formula, we find a natural isomorphism
\[\diag_*\so_X\FMcomp\s{E}={\R\pr_{1,3}}_*
(\pr_{2,3}^*\diag_*\so_X\lotimes\pr_{1,2}^*\s{E})\iso{\R\pr_{1,3}}_*
(\tilde{\diag}_*\so_{X\Timesa X}\lotimes\pr_{1,2}^*\s{E})\iso
{\R\pr_{1,3}}_*\tilde{\diag}_*\Ll\tilde{\diag}^*\pr_{1,2}^*\s{E}.
\]
On the other hand, $\so_{X\Timesa X}\FMcomp\s{E}\iso
{\R\pr_{1,3}}_*\pr_{1,2}^*\s{E}$, and it is clear that, with these
identifications, the morphism $\diag\mrs\FMcomp\id_{\s{E}}$
corresponds to
\[\R{\pr_{1,3}}_*\rladj{\tilde{\diag}}(\pr_{1,2}^*\s{E})\colon
{\R\pr_{1,3}}_*\pr_{1,2}^*\s{E}\to
{\R\pr_{1,3}}_*\tilde{\diag}_*\Ll\tilde{\diag}^*\pr_{1,2}^*\s{E}.\]
Now, ${\R\pr_{1,3}}_*\tilde{\diag}_*\Ll\tilde{\diag}^*\pr_{1,2}^*\s{E}
\iso\s{E}$ (because $\pr_{1,3}\comp\tilde{\diag}=
\pr_{1,2}\comp\tilde{\diag}=\id_{X\Timesa X}$), whereas
${\R\pr_{1,3}}_*\pr_{1,2}^*\s{E}\iso\pr_1^*\R{\pr_1}_*\s{E}$ (by the
``flat base change'' theorem for the Cartesian square \eqref{e33},
with $X$ in place of $Y$). To be more precise, the latter isomorphism
can be expressed as the composition
\[\pr_1^*\R{\pr_1}_*\s{E}\mor{\pr_1^*\R{\pr_1}_*\rladj{\pr_{1,2}}(\s{E})}
\pr_1^*\R{\pr_1}_*\R{\pr_{1,2}}_*\pr_{1,2}^*\s{E}\isomor
\pr_1^*\R{\pr_1}_*\R{\pr_{1,3}}_*\pr_{1,2}^*\s{E}
\mor{\lradj{\pr_1}(\R{\pr_{1,3}}_*\pr_{1,2}^*\s{E})}
\R{\pr_{1,3}}_*\pr_{1,2}^*\s{E},\]
where the middle map is the natural isomorphism (due to the fact that
$\pr_1\comp\pr_{1,2}=\pr_1\comp\pr_{1,3}$). Thus we can conclude that
$\Gi{1}\FMcomp\s{E}\iso
\cone{\R{\pr_{1,3}}_*\rladj{\tilde{\diag}}(\pr_{1,2}^*\s{E})}\iso
\cone{\lradj{\pr_1}(\s{E})}$, provided we show that in the diagram
(where the unnamed arrows denote the natural isomorphisms)
\[\xymatrix{\pr_1^*\R{\pr_1}_*\s{E}
\ar[rr]^-{\pr_1^*\R{\pr_1}_*\rladj{\pr_{1,2}}(\s{E})}
\ar[dd]_{\lradj{\pr_1}(\s{E})} 
\ar[drr] & & \pr_1^*\R{\pr_1}_*\R{\pr_{1,2}}_*\pr_{1,2}^*\s{E}
\ar[rr] & & \pr_1^*\R{\pr_1}_*\R{\pr_{1,3}}_*\pr_{1,2}^*\s{E}
\ar[d]^{\lradj{\pr_1}(\R{\pr_{1,3}}_*\pr_{1,2}^*\s{E})}
\ar[dll]_{\pr_1^*\R{\pr_1}_*\R{\pr_{1,3}}_*
\rladj{\tilde{\diag}}(\pr_{1,2}^*\s{E})\qquad} \\
& & \pr_1^*\R{\pr_1}_*\R{\pr_{1,3}}_*
\tilde{\diag}_*\Ll\tilde{\diag}^*\pr_{1,2}^*\s{E}
\ar[drr]_{\lradj{\pr_1}(\R{\pr_{1,3}}_*\tilde{\diag}_*
\Ll\tilde{\diag}^*\pr_{1,2}^*\s{E})\qquad}
& & \R{\pr_{1,3}}_*\pr_{1,2}^*\s{E}
\ar[d]^{\R{\pr_{1,3}}_*\rladj{\tilde{\diag}}(\pr_{1,2}^*\s{E})} \\
\s{E} & & & & 
\R{\pr_{1,3}}_*\tilde{\diag}_*\Ll\tilde{\diag}^*\pr_{1,2}^*\s{E}
\ar[llll]}\]
the outer square commutes. This follows from the fact that the three
inner triangles commute, as it can be easily checked using well known
compatibilities between adjunction morphisms.
\end{proof}

In order to apply the lemma for $\s{E}=\Gi{m}(m,0)$ we need to compute
$\R{\pr_1}_*\Gi{m}(m,0)$. {From} the definition \eqref{f1}
$\Gi{m}^0(m,0)\iso\diag_*\so_X(m)$, whereas for $j\ne0$
\begin{equation}\label{f3}
\Gi{m}^j(m,0)=\bigoplus_{-m<l\le0,\card{I}=-j-1,\w_I\le-l}
\so_{X\Timesa X}(l+m,-l-\w_I).
\end{equation}
To evaluate $\R{\pr_1}_*\Gi{m}(m,0)$ we use the following lemma, whose
proof is straightforward:

\begin{lemma}\label{f4}
With $X$ as above, and for two integers $p$ and $p'$, with $0<p<\sww$, we have
the isomorphism
\[
R^k{\pr_1}_*\so_{X\Timesa X}(p,p')\iso
H^k(X,\so_X(p))\otimes_{\K}\so_X(p')\iso
\begin{cases}
\p_p\otimes_{\K}\so_X(p') & \text{if $k=0$} \\
0 & \text{if $k\ne0$}.
\end{cases}
\]
\end{lemma}


Going back to \eqref{f3}, we see that $0<l+m\le m$. The inductive hypothesis
also assumed that $0<m<\sww$, thus $0<l+m<\sww$ and the lemma applies. Noting
that
\[
\R{\pr_1}_*\Gi{m}^0(m,0)\iso\R{\pr_1}_*\diag_*\so_X(m)\iso
{\pr_1}_*\diag_*\so_X(m)\iso\so_X(m),
\]
and applying Lemma~\ref{f4} we see that
every term of $\Gi{m}(m,0)$ is ${\pr_1}_*$-acyclic, hence
$\R{\pr_1}_*\Gi{m}(m,0)\iso{\pr_1}_*\Gi{m}(m,0)$. Inspecting the complex
${\pr_1}_*\Gi{m}(m,0)$ one observes that
\begin{equation}\label{f8}
\R{\pr_1}_*\Gi{m}(m,0)\iso{\pr_1}_*\Gi{m}(m,0)\iso\GimX{m},
\end{equation}
where for $0<m<\sww$ we defined
\begin{equation}\label{f5}
\Gim{m}:=\mc{{\pr_1}_*\pdr{1-m}(m,0)
\overset{\rho_m}{\too}\so_{\P}(m)},
\end{equation}
with each component of the morphism
\[
\rho_m^0 :\, {\pr_1}_*\pdr{1-m}^0(m,0)\iso
\bigoplus_{-m<l\le0}\p_{l+m}\otimes_{\K}\so_{\P}(-l)\too\so_{\P}(m)
\]
given as usual by the multiplication map.

Combining Lemma \ref{l1} and \eqref{f8}  we obtain that
\begin{equation}\label{f9}
\Gi{1}\FMcomp\Gi{m}(m,0)\iso
\mc{\mGi{m}\colon\so_X\boxtimes\GimX{m}\to\Gi{m}(m,0)},
\end{equation}
where $\mGi{m}^0\colon\so_{X\Timesa X}(0,m)\to\diag_*\so_X(m)$ is
the natural map, and for $j<0$
\begin{equation*}
\mGi{m}^j\colon
\bigoplus_{\substack{-m<l\le0\\ \card{I}=-j-1,\w_I\le-l}}
\p_{l+m}\otimes_{\K}\so_{X\Timesa X}(0,-l-\w_I)\too
\bigoplus_{\substack{-m<l\le0\\ \card{I}=-j-1,\w_I\le-l}}
\so_{X\Timesa X}(l+m,-l-\w_I)
\end{equation*}
is induced by the multiplication maps
$\p_{l+m}\otimes_{\K}\so_X\to\so_X(l+m)$.

It is easy to check that there is a natural morphism of complexes
$\mbo{m}\colon\bo{-m}\to\Gim{m}$, where $\mbo{m}^0=\id_{\so(m)}$ and
for $j<0$ each component of $\mbo{m}^j$
\[
\left(\mbo{m}^j\right)^I_{l,I'}\in
\Hom_{\P}(\so(m-\w_I),\p_{l+m}\otimes_{\K}\so(-l-\w_{I'}))
\iso\p_{l+m}\otimes_{\K}\p_{-m-l+\w_I-\w_{I'}}\]
(for $-m<l\le0$, $\card{I}=-j$, $\card{I'}=-j-1$, $\w_I\le m$ and
$\w_{I'}\le-l$) is given by
\begin{equation}\label{d2}
\left(\mbo{m}^j\right)^I_{l,I'}:=\begin{cases}
-(-1)^{N_I^i}\,x_i\otimes1\in\p_{\w_i}\otimes_{\K}\p_0 &
\text{if $l+m=\w_i$, $I=I'\cup\{i\}$} \\
0 & \text{otherwise}.
\end{cases}
\end{equation}

\begin{lemma}\label{Mres}
$\mbo{m}\colon\bo{-m}\to\Gim{m}$ is an isomorphism in $\D^b(\P)$ for any
$0<m<\sww$.
\end{lemma}

We will prove the lemma shortly, but first we look at its implications. Setting
\[
\amGi{m}=\mGi{m}\comp(\id_{\so_X}\boxtimes\,\mboX{m})\colon
\so_X\boxtimes\boX{-m}\too\Gi{m}(m,0),
\]
Lemma \ref{Mres} and \eqref{f9} imply that
\[
\Gi{1}\FMcomp\Gi{m}(m,0)
\iso\mc{\mGi{m}}\iso\mc{\amGi{m}} \in  \D^b(X\Times X).
\]

The components of $\amGi{m}$ can be calculated by composing
$\mGi{m}$ and $\id_{\so_X}\boxtimes\,\mboX{m}$.
After an explicit calculation it turns out that
$\amGi{m}^0\colon\so_{X\Timesa X}(0,m)\to\diag_*\so_X(m)$ is the
natural map; while for $j<0$
\[
\amGi{m}^j\colon (\so_X\boxtimes\boX{-m})^j\iso\so_X\boxtimes\boX{-m}^j\too
(\Gi{m}(m,0))^j\iso\pdrX{1-m}^{j+1}(m,0)
\]
can be identified with $\mpdrX{-m}^{j+1}(m,0)$. {From} the defining
equation \eqref{f1} of $\Gi{m}$ we see that
\[
\mc{\amGi{m}}^0\iso\Gi{m}^0(m,0)\iso\diag_*\so_X(m)\iso\Gi{m+1}^0(m,0),
\]
while for $j<0$, recalling \eqref{mcdiag},
\[
\mc{\amGi{m}}^j
\iso\Gi{m}^j(m,0)\oplus(\so_X\boxtimes\boX{-m})^{j+1}
\iso\pdrX{1-m}^{j+1}(m,0)\oplus\so_X\boxtimes\boX{-m}^{j+1}\iso
\pdrX{-m}^{j+1}(m,0)\iso\Gi{m+1}^j(m,0).
\]
Using the explicit form of $\amGi{m}$ it is also immediate to check that the
differential of $\mc{\amGi{m}}$ can be identified with that of $\Gi{m+1}(m,0)$.
Hence $\mc{\amGi{m}}\iso\Gi{m+1}(m,0)$ in $\CC^b(X\Times X)$, thereby proving
\eqref{e0}, and finishing the proof of the proposition.

\begin{proof}[Proof of Lemma \ref{Mres}]
We will show that $\cone{\mbo{m}}\iso0$, which immediately implies
that $\mbo{m}$ is an isomorphism in $\D^b(\P)$.
{From} the defining equation \eqref{f5} we have a distinguished triangle
\begin{equation*}\label{f6}
\xymatrix@1{
\Gim{m}[-1]\ar[r]^-{p} & {\pr_1}_*\pdr{1-m}(m,0) \ar[r]^-{\rho_m} &
\so(m) \ar[r] & \Gim{m},}
\end{equation*}
where $p$ is the natural projection. Using the octahedral axiom (TR4)
of triangulated categories we have
\begin{equation}\label{g1}
\begin{split}
\cone{\mbo{m}}&=\cone{\bo{-m}\mor{\mbo{m}}
\mc{{\pr_1}_*\pdr{1-m}(m,0)\mor{\rho_m}\so(m)}} \\
&\iso\cone{\cone{\bo{-m}[-1]\mor{p\comp\mbo{m}[-1]}
{\pr_1}_*\pdr{1-m}(m,0)}\too\so(m)}.
\end{split}
\end{equation}
Now, by \eqref{mcdiag}, for any $-\sww<k<0$ there is a distinguished triangle
\[\xymatrix@1{
\so(k)\boxtimes\bo{k}[-1]\ar[r]^-{\mpdr{k}} & \pdr{k+1} \ar[r] &
\pdr{k} \ar[r] & \so(k)\boxtimes\bo{k}}\]
in $\KK^b(\P\Times\P)$. Applying the exact functor
$\R{\pr_1}_*(\so(m,0)\otimes-)$ and reasoning as in the proof of 
Lemma~\ref{f4} gives another distinguished triangle in $\KK^b(\P)$:
\begin{equation*}\label{f7}
\p_{k+m}\otimes_{\K}\bo{k}[-1]\mor{{\pr_1}_*\mpdr{k}(m,0)}
{\pr_1}_*\pdr{k+1}(m,0)\too{\pr_1}_*{\pdr{k}}(m,0)
\too\p_{k+m}\otimes_{\K}\bo{k}.
\end{equation*}
{From} this we can deduce two things:
\begin{enumerate}
\item For $-\sww<k<-m$, since $\p_{k+m}=0$, it follows that
${\pr_1}_*\pdr{k+1}(m,0)\iso{\pr_1}_*\pdr{k}(m,0)$.
In other words
\[
{\pr_1}_*\pdr{-m}(m,0)\iso
{\pr_1}_*\pdr{-1-m}(m,0)\iso\cdots\iso{\pr_1}_*\pdr{1-\sww}(m,0).
\]
\item For $k=-m$ we have that
\[
\cone{{\pr_1}_*\mpdr{-m}(m,0)}\iso {\pr_1}_*\pdr{-m}(m,0).
\]
\end{enumerate}
{From} (1) and (2) it follows that in fact
\[
\cone{{\pr_1}_*\mpdr{-m}(m,0)}
\iso{\pr_1}_*\pdr{1-\sww}(m,0).
\]
On the other hand, it is easy to check that
$p\comp\mbo{m}[-1]\colon\bo{-m}[-1]\to{\pr_1}_*\pdr{1-m}(m,0)$ can
be identified with
${\pr_1}_*{\mpdr{-m}(m,0)}\colon\p_{0}\otimes_{\K}\bo{-m}[-1]\to
{\pr_1}_*\pdr{1-m}(m,0)$. Therefore
\[
\cone{p\comp\mbo{m}[-1]}\iso {\pr_1}_*\pdr{1-\sww}(m,0).
\]

An argument similar to the one used in the proof of
Lemma~\ref{f4} also shows that
$\FM[2]{\pdr{1-\sww}}(\so(m))\iso{\pr_1}_*\pdr{1-\sww}(m,0)$.
Combining this with the previous equation, and
using Prop.~\ref{diagres}, gives
\[
\cone{p\comp\mbo{m}[-1]}\iso\FM[2]{\pdr{1-\sww}}(\so(m))\iso
\FM[2]{\so_{\Delta}}(\so(m))\iso\so(m).
\]
By \eqref{g1} this means that $\cone{\mbo{m}}\iso
\cone{\so(m)\mor{\lambda}\so(m)}$ for some $\lambda\in\K$, and we claim that
$\lambda\ne0$. To prove this, observe that by \cite[Cor. 2.4.6]{Alberto} there
exists unique up to isomorphism $\s E\in\CC^b(\P)$ such that $\s
E\iso\cone{\mbo{m}}$ in $\D^b(\P)$, each $\s E^j$ is a sum of terms of the form
$\so(l)$ with $0\le l<\sww$ and $\s E$ is {\em minimal}, meaning that the
components of each $\diff{\s E}^j$ from $\so(l)$ to $\so(l)$ are $0$ for every
$l$. On the other hand, since $\mbo{m}^0=\id_{\so(m)}$ and for $j\ne0$ both
$\bo{-m}^j$ and $\Gim{m}^j$ are sums of terms of the form $\so(l)$ with $0\le
l<m$, it is easy to deduce that also $\s E^j$ for every $j\in\Z$ is a sum of
terms of the form $\so(l)$ with $0\le l<m$. So we must have $\lambda\ne0$
(otherwise obviously $\s E=\so(m)\oplus\so(m)[1]$), which proves that
$\cone{\mbo{m}}\iso0$.
\end{proof}

\begin{remark}\label{Fano}
{\em
Theorem~\ref{t1} admits the following generalization. Let $X$ be a smooth
stacky hypersurface of degree $d\le\w$ in $\P$. Then one can still define the
functors $\ms{L}$ and $\ms{K}$ as in \eqref{k2} and \eqref{k3}, although
$\ms{K}$ is not an equivalence when $d<\w$. Denoting by
$\langle\so_X(1),\dots,\so_X(\w-d)\rangle$ the full triangulated
subcategory of $\D^b(X)$ generated by the exceptional sequence
$(\so_X(1),\dots,\so_X(\w-d))$, we consider its right orthogonal
$\cat{D}:=\langle\so_X(1),\dots,\so_X(\w-d)\rangle^{\perp}$.
Explicitly, $\cat{D}$ is the full triangulated subcategory of
$\D^b(X)$ with objects the complexes $\s F$ such that
$\RHom_X(\so_X(i),\s F)=0$ for $0<i\le\w-d$ (obviously $\cat{D}=\D^b(X)$ if and
only if $d=\w$). It is easy to see that
$\ms{G}=\ms{L}\comp\ms{K}\colon\D^b(X)\to\D^b(X)$ sends $\cat{D}$ to itself,
and one proves that
\begin{equation}\label{eGr}
\underbrace{(\ms{G}\rest{\cat{D}}) \comp \cdots \comp
(\ms{G}\rest{\cat{D}})}_{\txt{$d$-times}} \iso (-)[2]
\end{equation}
(clearly this generalizes \eqref{eG}, which is the particular case $d=\w$).
When $\P=\P^{n+1}$ and $d<\sww=n+2$, \eqref{eGr} is equivalent to \cite[Lemma
4.2]{V14}, where $\ms{G}$ is replaced by $\ms{O}[1]=\ms{K}\comp\ms{L}$ and
$\cat{D}$ by $\cat{A}=\langle\so_X,\dots,\so_X(n+1-d)\rangle^{\perp}$ (the
equivalence of the two statements follows from the fact that
$\ms{G}\iso\ms{L}\comp\ms{O}[1]\comp\ms{L}^{-1}$ and
$\cat{D}=\ms{L}(\cat{A})$).

The proof of \eqref{eGr} is quite similar to that of Theorem~\ref{t1},
so here we just sketch the main steps, emphasizing what needs to be added
or changed. Of course, we still have $\ms{G}\iso\FM[1]{\G}$ with $\G$
defined as in Prop.~\ref{mainprop}, but now Claim~\ref{e3} holds only
for $0<m\le d$ (essentially with the same proof, and replacing
$\sww$ with $d$ in Lemma~\ref{f4}). On the other hand, Lemma \ref{l0}
must be changed to $\Gi{\sww}\iso\diag_*\so_X(-d)[2]$ (with the same
proof). Therefore \eqref{eGr} will follow if we prove that
\begin{equation}\label{dw}
\FM[1]{\Gi{d}(d,0)}\rest{\cat{D}}\iso
\FM[1]{\Gi{\sww}(d,0)}\rest{\cat{D}}.
\end{equation}
To this purpose observe that by what we have shown at the end of the
proof of Prop.~\ref{mainprop}, there are distinguished triangles in
$\KK^b(X\Times X)$ for $d\le m<\sww$
\[
\xymatrix@1{\so_X(d-m)\boxtimes\boX{-m} \ar[rr]^-{\amGi{m}(d-m,0)}
& & \Gi{m}(d,0) \ar[r] & \Gi{m+1}(d,0) \ar[r] &
\so_X(d-m)\boxtimes\boX{-m}[1]}.
\]
In order to conclude that \eqref{dw} holds, it is then enough to prove that
\begin{equation}\label{C0}
\FM[1]{\so_X(d-m)\boxtimes\boX{-m}}(\cat{D})=0
\end{equation}
for $d\le m<\sww$. Now, by part (1) of Lemma~\ref{kercomp}
\[
\FM[1]{\so_X(d-m)\boxtimes\boX{-m}}(\s F)\iso
\so_X(d-m)\otimes_{\K}\R\Gamma(X,\boX{-m}\otimes\s F)
\]
for every $\s F\in\D^b(X)$. Since $\bo{-m}\iso\abo{-m}$ in $\D^b(\P)$, where
each $\abo{-m}^j$ is a finite direct sum of terms of the form $\so_{\P}(i)$
with $m-\sww\le i<0$ (see \cite[Remarks 2.5.4 and 2.5.5]{Alberto}), we have
that
\[
\R\Gamma(X,\boX{-m}\otimes\s F)\iso\R\Gamma(X,\aboX{-m}\otimes\s F)\iso
\RHom_X(\aboX{-m}\dual,\s F),
\]
and the last term is $0$ if $d\le m<\sww$ (in which case clearly
$\aboX{-m}\dual\in\langle\so_X(1),\dots,\so_X(\w-d)\rangle$) and
$\s F\in\cat{D}$, thereby proving \eqref{C0}.

In analogy with \cite[Lemma 4.1]{V14} it can also be proved that the Serre
functor of $\cat{D}$ is isomorphic to $(\ms{G}\rest{\cat{D}})^{\comp
d-\sww}[n-1]$ (for $d=\sww$ this is just $(-)[n-1]$).
}
\end{remark}



\begin{thebibliography}{10}

\bibitem[A]{Paul:TASI2003}
P.~S. Aspinwall,
\newblock {\em D-branes on Calabi-Yau manifolds},
\newblock in ``Recent Trends in String Theory'', pages 1--152, World
  Scientific, 2004,
\newblock hep-th/0403166.

\bibitem[AHK]{en:Horja}
P.~S. Aspinwall, R.~P. Horja, and R.~L. Karp,
\newblock {\em Massless D-branes on Calabi-Yau threefolds and monodromy},
\newblock Commun. Math. Phys. {\bf 259} (2005) 45--69, hep-th/0209161.

\bibitem[AKO]{Auroux}
D.~Auroux, L.~Katzarkov, and D.~Orlov,
\newblock {\em Mirror symmetry for weighted projective planes and their
  noncommutative deformations},
\newblock math.AG/0404281.

\bibitem[B]{Bei:res}
A.~A. Beilinson,
\newblock {\em Coherent sheaves on $\P^n$ and problems in linear algebra},
\newblock Funktsional. Anal. i Prilozhen. {\bf 12} (No. 3)  (1978) 68--69.

\bibitem[BCS]{BorisovDM}
L.~A. Borisov, L.~Chen, and G.~G. Smith,
\newblock {\em The orbifold {C}how ring of toric {D}eligne-{M}umford stacks},
\newblock J. Amer. Math. Soc. {\bf 18} (No. 1)  (2005) 193--215,
  math.AG/0309229.

\bibitem[Br1]{Bridgeland:2002}
T.~Bridgeland,
\newblock {\em Stability conditions on triangulated categories},
\newblock math.AG/0212237.

\bibitem[Br2]{Bridgeland:Review}
T.~Bridgeland,
\newblock {\em Derived categories of coherent sheaves},
\newblock math.AG/0602129.

\bibitem[C]{Alberto}
A.~Canonaco,
\newblock {\em The Beilinson complex and canonical rings of irregular
  surfaces},
\newblock Mem. Amer. Math. Soc. {\bf 183} (No. 862)  (2006),
math.AG/0610731.

\bibitem[CK]{Cox:Katz}
D.~A. Cox and S.~Katz,
\newblock {\em Mirror symmetry and algebraic geometry}, Mathematical Surveys
  and Monographs~{\bf 68},
\newblock AMS, Providence, RI, 1999.

\bibitem[Clay]{DBook}
P.~S. Aspinwall et al.,
\newblock {\em Dirichlet branes and mirror symmetry},
\newblock Clay mathematics monographs (2008).

\bibitem[D]{Douglas:2000gi}
M.~R. Douglas,
\newblock {\em D-branes, categories and N = 1 supersymmetry},
\newblock J. Math. Phys. {\bf 42} (2001) 2818--2843, hep-th/0011017.

\bibitem[H]{Hart:dC}
R.~Hartshorne,
\newblock {\em Residues and duality},
\newblock Lecture notes of a seminar on the work of A. Grothendieck, given at
  Harvard 1963/64. Lecture Notes in Mathematics, No. 20, Springer-Verlag,
  Berlin, 1966.

\bibitem[Ho1]{Horj:DX}
R.~P. Horja,
\newblock {\em Hypergeometric Functions and Mirror Symmetry in Toric Varieties}
  (1999),
\newblock math.AG/9912109.

\bibitem[Ho2]{Horj:EZ}
R.~P. Horja,
\newblock {\em Derived category automorphisms from mirror symmetry},
\newblock Duke Math. J. {\bf 127} (No. 1)  (2005) 1--34,
  math\-.\-AG\-/\-0103231.

\bibitem[Hu]{Huybrechts}
D.~Huybrechts,
\newblock {\em {F}ourier-{M}ukai Transforms in Algebraic Geometry},
\newblock Oxford Mathematical Monographs, 2006.

\bibitem[K1]{K1}
R.~L. Karp,
\newblock {\em ${\mathbb C}^2/{\mathbb Z}_n$ fractional branes and monodromy},
\newblock Commun. Math. Phys. {\bf 270} (2007) 163--196, hep-th/0510047.

\bibitem[K2]{K2}
R.~L. Karp,
\newblock {\em On the ${\mathbb C}^n/{\mathbb Z}_m$ fractional branes},
\newblock hep-th/0602165.

\bibitem[Ka]{Kawamata:DC}
Y.~Kawamata,
\newblock {\em Equivalences of derived categories of sheaves on smooth stacks},
\newblock Amer. J. Math. {\bf 126} (No. 5)  (2004) 1057--1083, math.AG/0210439.

\bibitem[Ki]{King}
A. King,
\newblock {\em Tilting bundles on some rational surfaces},
\newblock 1997, http://www.maths.bath.ac.uk/~masadk/papers/tilt.ps,
unpublished.

\bibitem[Ko]{Kont:mon}
M.~Kontsevich,
\newblock {\em Lecture at Rutgers University},
\newblock unpublished (Nov. 11, 1996).

\bibitem[Ku]{V14}
A.~Kuznetsov,
\newblock {\em Derived categories of cubic and $V_{14}$ threefolds},
\newblock math.AG/0303037.

\bibitem[LM]{Stacks}
G.~Laumon and L.~Moret-Bailly,
\newblock {\em Champs alg\'ebriques}, Ergebnisse der Mathematik und ihrer
  Grenzgebiete.~{\bf 39},
\newblock Springer-Verlag, Berlin, 2000.

\bibitem[M]{Mor:geom2}
D.~R. Morrison,
\newblock {\em Geometric aspects of mirror symmetry},
\newblock in ``Mathematics unlimited---2001 and beyond'', pages 899--918,
  Springer, Berlin, 2001,
\newblock math.AG/0007090.

\bibitem[ST]{ST:braid}
P.~Seidel and R.~Thomas,
\newblock {\em Braid group actions on derived categories of coherent sheaves},
\newblock Duke Math. J. {\bf 108} (No. 1)  (2001) 37--108, math.AG/0001043.

\end{thebibliography}
\end{document}